# Best Practices on QSP Model Reporting for Regulatory Use: perspectives from ISoP QSP SIG Working Group

Susana Zaph* (Sanofi), Blerta Shtylla* (Pfizer), Steve Chang (Sims+), Yougan Cheng (Daiichi Sankyo), Jingqi Q.X. Gong (GSK), Abhishek Gulati (Merck & Co., Inc.), Emma Hansson (Pharmetheus), Alexander Kulesza (ESQlabs), Alexander V. Ratushny (Bristol Myers Squibb), Federico Reali (COSBI), Conner Sandefur (Sims+), Brian Schmidt (Madrigal), Fulya Akpinar Singh (Genmab), Monica Susilo (Genentech), Weirong Wang (JnJ)

*Corresponding Author(s)

Abstract

Quantitative systems pharmacology (QSP) models are increasingly applied to inform decision-making across drug development and to support regulatory interactions within model-informed drug development (MIDD). QSP supports a broad range of applications across drug development and can be tailored to specific therapeutic areas, mechanisms of action, and contexts of use (CoU). While this diversity is a core strength of QSP, it also presents challenges for reporting for regulatory use. Despite the growing impact of QSP models, there is currently no established guidance on how QSP analyses should be documented and reported for regulatory purposes. This white paper, developed by the International Society of Pharmacometrics (ISoP) QSP Special Interest Group Working Group on Credibility Assessment of QSP for Regulatory Use, seeks to address this gap by proposing best practices for QSP model reporting in regulatory settings. The recommendations are grounded in collective real-world experience from regulatory interactions and are aligned with reporting guidance established for physiologically based pharmacokinetic (PBPK) modeling and reporting principles outlined in ICH M15. Rather than prescribing a rigid, one-size-fits-all template, this work proposes a flexible, tiered reporting framework that accounts for development phase and model impact. The proposed framework is intended to facilitate regulatory review and enhance transparency while accommodating the inherent diversity of QSP modeling.

## 1 Introduction to QSP Modeling in Regulatory Context

Quantitative Systems Pharmacology (QSP) modeling has emerged as a useful tool in model-informed drug development (MIDD), by integrating pre-existent knowledge to create a representation of the disease and drug response informed by diverse data sources. This modeling approach can offer mechanistic insights that complement more empirical models and support a variety of applications both for internal pharma or biotech decision-making and more recently, regulatory use. In early discovery, QSP assists in prioritizing targets and predicting efficacy; during preclinical and translational phases, it aids in dose selection; and throughout clinical development, it informs trial design, patient stratification, and regulatory interactions.

QSP approaches are typically tailored to specific therapeutic areas, mechanisms of action, and stages of drug development, leading to a significant amount of model diversity. These mechanistic models can range from simpler compartmental systems to complex multi-scale representations involving sub-cellular and cellular signaling, tissue-level dynamics, and even whole-body physiology. QSP models can be built using diverse mathematical approaches such as differential equations, agent-based simulations, or hybrid frameworks, depending on the biological questions and data availability. This adaptability enables QSP modeling to address emerging drug development challenges such as combination therapies and variability in patient response. This versatility is a core strength, allowing QSP to be tailored to specific scientific questions.

The diversity of QSP approaches and applications is a major strength, but it also creates challenges to standardization. QSP models are often built with a specific context of use (CoU) in mind, incorporating only the biological mechanisms and data relevant to that purpose. As a result, there is no universal framework or modeling template that fits all use cases. For example, a model designed for immuno-oncology may involve complex immune cell interactions driven by receptor engagement with tumor cells, whereas a model for metabolic disease might focus on hormonal regulation and organ-level feedback loops. These differences lead to variability in model structure, mathematical formalisms, and data requirements for calibration and validation, often extending beyond the data routinely collected during a clinical trial. Thus, the same diversity that underpins QSP's value also complicates the establishment of strict guidelines for model qualification, reproducibility, and reporting.

Regulatory agencies increasingly recognize QSP as a component of model-informed drug development (MIDD) (1), with growing use in regulatory decision-making. Two FDA landscape analyses (2) (3) have described the expanding therapeutic applications and diverse contexts of use across different types of submissions including NDAs and BLAs.

In the absence of regulatory QSP reporting guidelines, the QSP community and regulatory agencies have often relied on related mechanistic modeling approaches such as physiologically based pharmacokinetic (PBPK) for model assessment (4) as well as reporting frameworks (5) (6). Recent finalization of the International Council for Harmonization of Technical Requirements for Pharmaceuticals for Human Use (ICH) M15 harmonized guideline (1) establishes high-level expectations for MIDD planning, evaluation, and documentation, including principles applicable to mechanistic models. However, it does not provide detailed guidance tailored to specific considerations for QSP.

Because QSP models are mechanistically complex and structurally diverse, more specific reporting standards can be beneficial. Clearly communicating the validity of mathematical representations, embedded biological assumptions within each CoU remain challenging, particularly for multi-scale differential equation systems implemented in heterogeneous modeling platforms (e.g., MATLAB, R, Julia, SimBiology). Unlike PBPK modeling, which often relies on more standardized tools such as Simcyp (7), G+, PKSim, QSP lacks computational platform harmonization, complicating model sharing, reproducibility and regulatory review. Given tight submission timelines, ensuring transparent justification of model structure, parameters, uncertainty, and validation is essential to support regulatory confidence and maximize the impact of QSP analyses.

To address these challenges, the QSP community would benefit from gathering and developing assessment and reporting best practices from the growing number of regulatory QSP use cases in recent years. Future reporting framework elements should be flexible and proportionate, enabling efficient evaluation based on CoU and model impact. Given the diversity and complexity of QSP models, reporting components may require careful planning and would benefit if they are developed and integrated organically, with increasing depth aligned with drug development stage. Such a reporting framework, when used routinely, could enhance multi-disciplinary team communication and collaboration while supporting regulatory review by promoting consistent level of documentation appropriately tailored to regulatory milestone and reducing technical barriers to the interpretation of QSP analyses.

In 2024, a Working Group focused on "Credibility Assessment of QSP for Regulatory Use" was launched under the umbrella of the International Society of Pharmacometrics (IsoP) QSP Special Interest Group (SIG) with representatives from QSP teams embedded in diverse organizations such as pharmaceutical and biotechnology companies, contract research organizations and academic institutions. The working group collected and synthesized best practices from their respective organizations as well as published case studies.

To build a QSP reporting best practices framework, a two-step approach was applied: 1) review and evaluation of special considerations that might be applicable to QSP from existing relevant regulatory reporting guidance (e.g. PBPK, ICH M15) and 2) synthesis of existing reporting practices from collected case studies of health authority interactions that utilized QSP models. This white paper summarizes QSP reporting best practices from the Working Group that can facilitate the use of QSP in various regulatory interactions. This paper draws on case studies of QSP applications that have informed regulatory decisions across diverse therapeutic areas and CoU, reflecting the breadth of QSP models. It also incorporates reporting elements expected to gain importance as the field matures and relevant guidelines, such as ICH M15, are implemented. By improving the clarity and consistency of QSP model reporting in regulatory submissions, the Working Group hopes that these best practices will help transition QSP from a specialized modeling approach to a mainstream tool of model informed drug development.

# 2 Considerations for QSP reporting elements from ICH M15 and PBPK reporting guidelines

The ICH M15 guideline (1) introduced a structured framework for assessing and reporting MIDD evidence to inform regulatory decision-making. Six key model assessment elements were introduced to enable a structured, risk-based assessment of MIDD evidence that enables transparency and clear understanding of the proposed MIDD strategy and implementation: 1) the question of interest, 2) CoU, 3) model influence, 4) consequence of wrong decision, 5) model risk, and 6) model impact. In addition to the key assessment elements, additional considerations for interaction with regulators are specified to support planning and regulatory interactions: technical criteria, appropriateness of proposed MIDD, evaluation of the model(s) and outcomes, and the outcome of MIDD evidence assessment. The framework will support development of best practices for QSP model assessment, facilitate transparent and efficient dialogue with regulatory authorities, and enable timely and consistent communication between QSP modelers and multidisciplinary teams on the model's intended scope and impact on drug development. Finally, the ICH M15 guideline proposed a structure for Model Analysis Reports (MARs) which the Working Group believes establishes an overarching framework that can be customized to meet the specific needs for QSP reporting.

The FDA (5) and EMA (6) have also published PBPK guidelines earlier and despite differences in terminology and format, the agencies seem aligned on the foundational principles for reporting of PBPK. The Working Group evaluated whether components from these guidelines could be suitable to support specific considerations for QSP model reporting within the ICH M15 MAR structure. Specifically, even though the current PBPK guidelines do not explicitly use the ICH M15 terminology for MIDD assessment, some of the core concepts are embedded within them and can be relevant in contextualizing QSP assessment. Both agencies require a clear statement of the analysis objectives and regulatory rationale (e.g., informing dose selection, supporting labeling, replacing or augmenting clinical data, evaluating DDIs, or assessing special populations). However, the treatment of model risk differs in specificity between agencies. The EMA guidance explicitly addresses model qualification requirements based on regulatory impact, which states that higher regulatory impact (e.g., using PBPK model in place of clinical data) requires stricter model qualification. In contrast, the FDA PBPK reporting guidance refers to a risk-based approach for credibility assessment aligned with the ASME V&V40 and exemplified by Kuemmel et. al. (4) for PBPK. Moving forward, best practices for QSP reporting will incorporate the six assessment elements defined in the ICH M15 to ensure a consistent framework for evaluating model risk and impact, and these elements should be clearly listed in a table as the MIDD evidence assessment as recommended in the ICH M15. Sharing this table as an appendix to the MAR is recommended.

**Table 1** is a summary of key reporting elements in the PBPK guidelines and the ICH M15 guidance and the Working Group's perspective of how these principles may apply to QSP reporting. Section 4 describes in more detail the QSP considerations of each reporting element.

TABLE 1 MODEL ANALYSIS REPORT COMPONENTS CONSIDERATIONS BASED ON PBPK REPORTING AND ICH M15 GUIDANCE CONTENT

| Sections | Content based on PBPK and ICH M15 | QSP Considerations |
|---|---|---|
| **Executive Summary** | • An overview of the rationale for the analyses<br>• A brief summary of data and methods<br>• A brief summary of the results and conclusions | • *No specific QSP considerations* |
| **Introduction** | • The rationale for the analyses<br>• Relevant ***background information and knowledge**** <br>• Description of prior analyses with reference to prior reports<br>________<br>** The FDA and EMA PBPK guidelines recommend providing background information describing the molecule and the modeling context. This includes the drug's physicochemical characteristics, PK/PD properties, and a summary of the exposure-response (ER) relationship for efficacy and safety. However, the agencies emphasize different aspects of this background. The EMA places a stronger emphasis on technical depth, specifically requesting a quantitative mass-balance diagram and extensive background information on the model itself. In contrast, the FDA emphasizes regulatory context, prioritizing prior PBPK-related regulatory interactions and cross-referencing previous PBPK study reports.* | • The rationale for the analyses<br>• Relevant ***background information and knowledge**** <br>• Description of ***prior analyses with reference to prior reports***** <br>________<br>** For QSP reporting, the emphasis could shift toward clearly describing the formalized biological mechanisms and their effect on the pathophysiology of interest. This may include outlining the biological knowledge and evidence represented in the model, the rationale for the selected model scope and architecture, and how CoU or data availability influences that scope. Analogous to a mass-balance diagram, a schematic capturing key biological interactions and illustrating how the drug's MoA mechanistically explains the exposure-response relationship would be helpful.*<br>*** It may also be advisable to publish an initial version of a novel mathematical formulation developed to represent a new MoA and reference it accordingly.*<br>*See Section* 4.2.1 *for more details* |
| **Objectives** | • The objectives of the analysis include the intended application of the model | • *No specific QSP considerations* |
| **Data and Methods** | *Data Sources*<br>• Criteria and rationale wrt source data inclusion/exclusion<br>• Relevant design features of studies and/or experiments<br>*Model and Simulation Details:*<br>• Model development, assumptions†<br>• Methods, computational platforms‡ | *Data Sources*<br>• *No specific QSP considerations*<br>*Model and Simulation Details:*<br>• ***Model development, assumptions†***<br>• ***Methods, computational platforms‡*** |
| | • Strategic approaches i.e. sequence of development, numerical methods⁂<br><br>*Approaches for Model Evaluation:*<br>• Verification[a], validation[b] and applicability assessment[c]<br>• Uncertainty assessment is recommended<br>• Prediction and simulation methods and scenarios (if relevant)<br>• Detailed technical criteria for model evaluation and model outcomes<br><br>———<br>*†<u>Model Parameters</u>: Both PBPK guidelines require listing all model parameters along with their sources or scientific justification.*<br>*‡<u>Model Simulation</u>: Both guidelines require the description of the simulation design for the intended scenarios.*<br>*⁂ <u>Model Workflow</u>: Both PBPK guidelines require a description of the model development steps, including model construction, parameter loading, verification, modification, and evaluation prior to application. The EMA explicitly provides an example of a model workflow diagram in their guidance. The FDA suggests outlining model development procedures using a workflow, decision tree, table, or other schematic representation. Although current PBPK guidance from both agencies do not specifically use the same terminology, they listed similar concepts. These reporting elements should be considered in QSP reports.*<br>*a. <u>Model Verification:</u> Verification step ensures that the model implementation is accurate.*<br>*b. <u>Model Validation:</u> The validation step focuses on assessing the robustness of the model's predictive performance against real-world observations. Validation strategy should be tailored to the CoU. The EMA guidance recommends presenting comparisons using goodness of fit plots (plots of simulated vs. observed data) along with tables and visual predictive plots (VPCs). The EMA specifically emphasizes that acceptance criteria (adequacy of prediction) for the agreement of simulated with observed data depend on* the *regulatory* | • ***Strategic approaches i.e. sequence of development, calibration strategy numerical methods⁂***<br>*Approaches for Model Evaluation:*<br>• ***Verification[a], validation[b] and applicability assessment[c]***<br>• Uncertainty assessment is recommended<br>• Prediction and simulation methods and scenarios (if relevant)<br>• Detailed technical criteria for model evaluation and outcomes<br><br>———<br>*†<u>Model Parameters</u>: For QSP models, providing a comprehensive table that includes parameter values, units, data sources, and their biological interpretation is essential for robust model evaluation. QSP model parameters may be derived from diverse data types (in vitro, in vivo, ex vivo, clinical, or inferred) and their provenance should be considered when assessing their biological relevance and plausibility.*<br>*‡ <u>Model Simulation</u>: What simulations were performed, the rationale for selecting specific scenarios or virtual populations, and how simulation outputs inform the model's CoU. Assumptions and key inputs should be clearly specified for each scenario.*<br>*⁂<u>Model Workflow</u>: Recommended workflow elements include calibration and validation strategies, mapping of datasets used to inform model components, and approaches for virtual population generation.*<br>*a. <u>Model Verification</u>: Code verification of correct implementation of equations, units, and parameter definitions. Structural verification evaluates whether mathematically represented biological processes behave as intended, confirming the absence of issues such as mass-balance, violations or numerical artifacts.*<br>*b. <u>Model Validation:</u> QSP validation often spans multiple biomarkers, pathway dynamics, clinical endpoints, frequently requiring VPop-based comparisons rather than a single "typical" trajectory. Acceptance criteria should consider biological heterogeneity and measurement uncertainty and should specify appropriate observed dataset comparisons and metrics.* |
| | *impact and need to be considered separately for each CoU. The consequences of poor prediction should be discussed and justified, e.g., the acceptance limits for a victim drug must be set in the context of the concentration-effect and concentration-safety relationships of the drug at the time of submission. Biologically plausible explanations for any discrepancy between predicted and observed data should also be considered. Sensitivity analysis was mentioned in both the FDA and EMA guidances.*<br>*c. Applicability Assessment: The assessment focuses on the adequacy of the data and model for the intended use of the model.* | *Validation should be justified based on MIDD evidence assessment and aligned with model influence and decision consequence.*<br>*c. Applicability Assessment: Applicability for QSP models depends on whether the represented biology is appropriate for the intended CoU, whether parameter uncertainty is acceptable for the decision context, and whether VPop construction is adequately justified to support modeling conclusions.*<br><br>*See Section 4.2.2 for more details* |
| **Results** | • Results of model development and model evaluation, including predictions and simulations (e.g., parameter estimates and associated uncertainty) as graphical and/or tabular displays<br>• Detailed results of the assessment against the technical criteria for model evaluation and interpretation of model outcomes<br>• Where relevant, deviations from the Model Analysis Plan should be described and justified. | • Results of model development and model evaluation, including predictions and simulations (e.g., parameter estimates and associated uncertainty) as graphical and/or tabular displays<br>• Detailed results of the assessment against the technical criteria for model evaluation and interpretation of model outcomes<br>• ***Where relevant, deviations from the Model Analysis Plan should be described and justified*** [d]<br><br>d. *Model analyses plans are not routinely used for QSP but are supported as part of early planning and regulatory engagement. Peer review publication of models especially describing first of its kind formulation for novel modalities or less studied pathophysiology is recommended.* |
| **Discussion** | • Interpretation of results, including data and model adequacy, limitations of the data and model, and clinical and/or other implications | • *No specific QSP considerations* |
| **Conclusions** | • The conclusions of the analyses | • *No specific QSP considerations* |
| **Appendixes** | Additional materials cross-referenced in the MAR may include, for example:<br>• A reference list covering data sources used in the analyses (e.g., bioanalytical reports, clinical study reports, laboratory reports, or literature)<br>• Supplemental data descriptions and model development and evaluation | Additional materials cross-referenced in the MAR may include, for example:<br>• A reference list covering data sources used in the analyses (e.g., bioanalytical reports, clinical study reports, laboratory reports, or ***literature***)<br>• Supplemental data descriptions and model development and evaluation |
| | results, including graphical and/or tabular displays, as appropriate<br>• User-generated code for the relevant model(s)[e]<br><br>e. *Both the FDA and EMA require the electronic submission of all relevant model files in an executable format to ensure reproducibility. The EMA guideline recommends that the files be listed in a tabular format within the study report. For FDA guidance, files must adhere to current eCTD specifications regarding file types and module locations. The FDA also encourages the inclusion of an orientation document to guide reviewers through the files, along with supporting data (e.g., clinical PK/PD) and clear, hyperlinked cross-references to other relevant sections of the regulatory dossier (e.g., IND, NDA).* | results, including graphical and/or tabular displays, as appropriate<br>• ***User-generated code for the relevant model(s)***[e]<br><br>e. *User Generated Code for QSP: The report should specify the software used for QSP with analyses, including exact versions and all required libraries or packages to ensure reproducibility. To facilitate reviewer understanding, an orientation document may include a workflow diagram illustrating relationships among datasets, scripts, and analyses, along with a step-by-step guide for running simulation workflows. Scripts should be provided in an executable format with clear instructions for installation of required packages, libraries or toolboxes. Consistent with best practices, scripts should be modular, well-commented, and version-controlled, with each script corresponding to a specific analysis or figure included in the QSP report*. |

# 3 Overview of documentation and reporting components from a survey and selected case studies

QSP model scope can span multiple biological scales and incorporate varying levels of mathematical complexity. Model scope is typically refined as programs transition across milestones from discovery to late development. During this time, which can span multiple years, the non-clinical and published datasets supporting key model assumptions often evolve. Considering this, the Working Group recommends incremental and consistent development of model documentation (including curation of user generated codes) that can in parallel support reporting needs for a variety of regulatory interactions. Given that QSP models integrate novel biological components and can rely on novel mathematical formulations it is beneficial for teams to plan early for publication and external peer review. This provides early vetting of key assumptions by domain area expert peers, strengthening scientific credibility. Publication as a step in building QSP model documentation can facilitate robust and streamlined MAR generation and facilitate transparent communication with regulatory reviewers, illustrated in Figure 1.

QSP models and corresponding analyses may support various regulatory milestones, each with its own reporting needs: from early interactions to all the way into submission and post-approval activities, Figure 1. The landscape FDA reviews (2) (3) have been helpful to the QSP community since they illustrate the increasing number of FDA submissions, with most leveraging QSP to address efficacy (66.2%) and safety (33.8%) and a strong emphasis on dose optimization strategies. In addition, submissions in multiple therapeutic areas have included QSP analyses with oncology leading QSP adoption (accounting for 50% of all FDA submissions). A variety of CoU were reported including dose optimization across clinical development, pediatric extrapolation, safety risk mitigation, NDA/BLA supplements for new indications and dosing regimens. The reported increased integration of QSP into regulatory submission reflects its maturation as a decision-supporting tool but also highlights the need for the key stakeholders in industry, academia and regulatory authorities to collaborate and develop best practices on MIDD evidence assessment and reporting of QSP; these conversations are facilitated tremendously by sharing of case studies that highlight current practices for QSP reporting and feedback from regulators based on the nature of their review.

FIGURE 1 OVERVIEW OF RESEARCH, DEVELOPMENT AND REGULATORY INTERACTION MILESTONES WITH EXISTING DOCUMENTATION PRACTICES FOR QSP

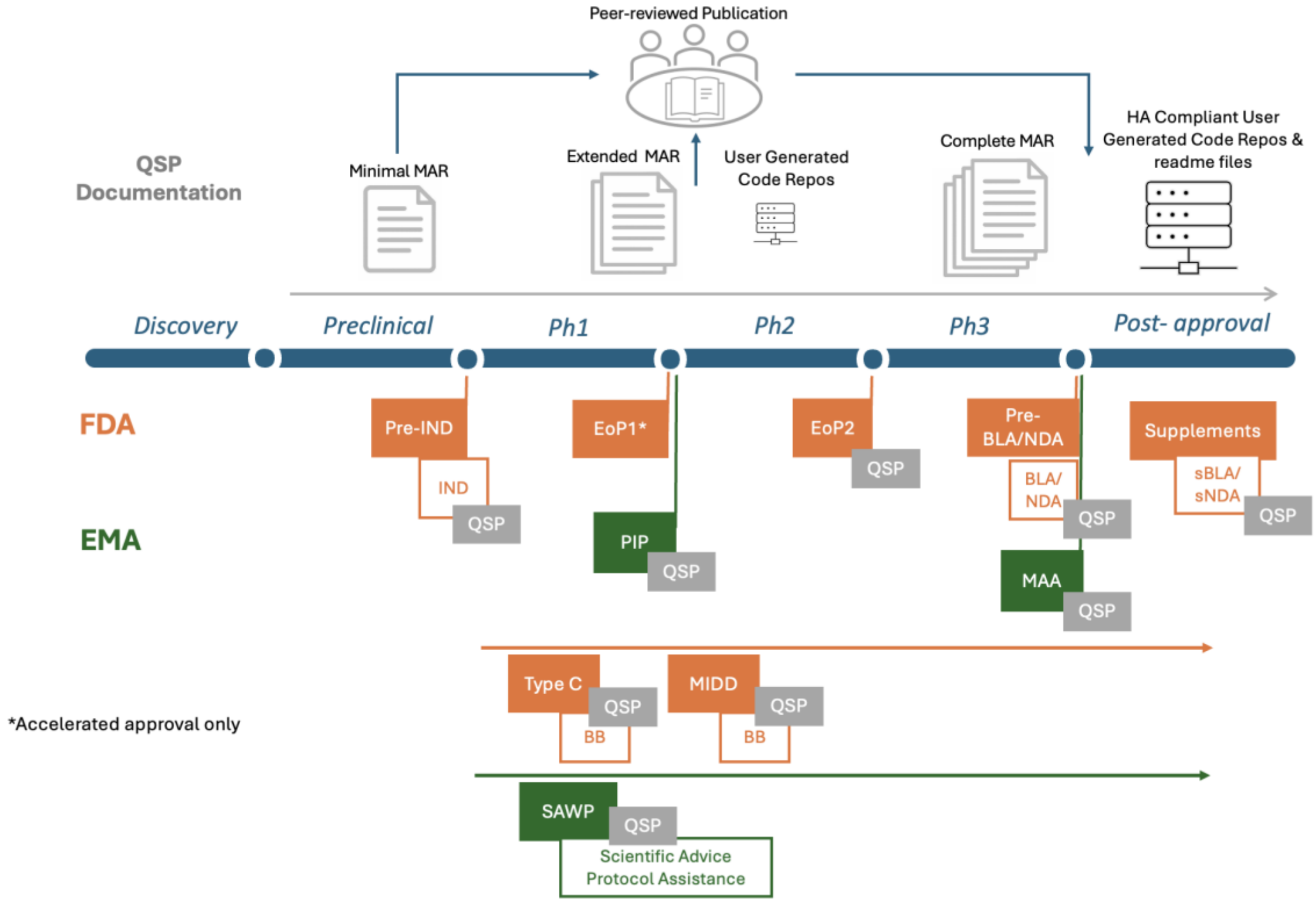

In the earlier stages of drug development, QSP regulatory use has expanded rapidly. Unlike BLA/NDA submission stage, in early clinical development data may be limited and the characterization of the drug's exposure-response relationships for efficacy and safety may still be emerging, so documentation practices can vary. In this setting, examples of QSP CoU could be first in human (FIH) dose projections, recommended phase 2 or 3 dose (RP2D/RP3D) optimization/justification; mechanistic biomarker-based patient stratification or providing mechanistic insight into combination treatment or dose optimization across indications (8) (9).

## 3.1 Summary of regulatory QSP reporting elements collected from Working Group

The experiences from the Working Group with regards to reporting practices of QSP used in regulatory interactions across the research and development spectrum were collected via an anonymized informal survey and summarized in Table 2. The intent of this survey was to identify common components included in QSP reporting and the frequency in which these components appear in various regulatory interactions for Working Group members.

TABLE 2 SUMMARY OF TYPES AND FREQUENCY OF QSP REPORTING ELEMENTS BASED ON WORKING GROUP EXPERIENCES

| QSP reporting format | Non BLA/NDA regulatory interaction | BLA/NDA submission | MIDD meeting |
|---|---|---|---|
| Concise summary in protocol or briefing book or module 2.7.2 | Always | Always | Always |
| Minimal MAR as appendix | Sometimes | Rarely | Rarely |
| Complete MAR as appendix | Rarely | Usually | Usually |
| Risk-based credibility assessment | Always | Always | Always |
| Model code repository and readme files | Never | Always | Usually |

Always Usually Sometimes Rarely Never

In addition, a set of anonymized working group case studies as well as FDA published multidisciplinary reviews of submissions with QSP analyses was collected and evaluated to help identify existing practices for QSP reporting across development stages and type of regulatory interactions. A *case study evaluation rubric* was developed to facilitate consistent evaluation of case studies and reduce potential subjectivity in evaluation criteria. The rubric is included in Appendix Figure 1 and included questions to help identify the presence of reporting elements that could describe the following five categories: 1) model purpose and regulatory context, 2) model structure and complexity, 3) model development and validation, 4) model results and conclusions, 5) regulatory feedback and impact. Each case study was independently reviewed by a working group member, and then rubric notes were shared and discussed with the rest of the group to help fill gaps and clarify interpretation as well as complement information gaps with publications from case study models that had been publicly disclosed. In view of the utility of publishing in the process of building QSP model documentation and as well as helping provide details on model development, a specific question with regards to status of publication of the QSP model was included in the rubric's first category.

### 3.1.1 Non BLA/NDA regulatory interaction reporting elements

Based on the informal survey, the CoU and the degree of model influence were the primary determinants of the level of detail required for QSP reporting in early-phase regulatory communications. The extent and format of QSP reporting in IND submissions depended on the specific CoU being addressed and on how central the QSP analysis was to inform dose selection. Across all IND-related working group survey responses, a brief description of QSP results was consistently included in the section describing the methods used to determine dose selection in the study protocol.

Whether a MAR was required depended on several factors. When relevant clinical data are limited—such as during first-in-human (FIH) studies or in novel drug combination programs supported only by single-agent clinical data—QSP analyses can provide essential mechanistic insight to address key translational questions. As the influence of the QSP model increases under these circumstances, appropriate validation of the model using comparator clinical data most relevant to the CoU becomes increasingly important. Conversely, when QSP analyses serve as supportive evidence to complement existing clinical data (e.g., for RP2D justification), reporting can be more concise. In such cases, reporting may include a high-level summary of the QSP analysis, referred to here as a *minimal MAR*, provided as an appendix. A *complete MAR*, which includes more comprehensive model documentation, is typically reserved for select cases—particularly when the QSP analysis serves as the primary basis for dose selection or when the CoU is atypical (e.g., back-translation of efficacy to a related indication). The reporting elements that define a minimal MAR vs a complete MAR are listed in Section 4.

Two anonymized early-development case studies were evaluated: one representing an IND submission for a biologic, and the other supporting a regulatory interaction for dose selection for a Phase 2 proof-of-concept (POC) study (Appendix Figure 2). Across these examples, an increase in reporting detail was observed with advancing drug development stage. The IND submission included a minimal MAR with a brief description of model development and data sources, whereas the Phase 2 case study MAR provided more detailed descriptions of the calibration and validation strategies. Similarly, model limitations were clearly stated in both the IND and Phase 2 QSP MAR. Local sensitivity analyses and virtual population approaches were briefly described in both cases to support exploration of parameter sensitivity and output heterogeneity, highlighting that these elements constitute minimal components that should be included in all types of QSP reports. For both early submissions, the base models had not yet been published; however, for the Phase 2 case study, a model publication was planned, consistent with the Working Group's recommendation to seek peer review for model components that can be disclosed using published datasets.

Taken together, these case studies support the recommendation that core documentation elements be included in all QSP reports, independent of the type of regulatory interaction. These elements include model description and limitations, data sources, sensitivity analysis and/or virtual population sampling, and risk matrix assessment. Based on the collective experience of the Working Group, when MARs contain an adequate level of detail and are deemed low risk, inclusion of user-generated code is generally not necessary. However, if requested by reviewers, QSP teams should be prepared to disclose code able to sufficiently reproduce MAR components to facilitate efficient regulatory review. Overall, these case studies illustrate a shared best-practice experience: documentation and reporting depth naturally increase as programs mature and as additional data and supporting evidence become available.

### 3.1.2 FDA MIDD or EMA scientific advice meeting reporting elements

Special considerations apply for the FDA MIDD paired meetings or EMA's scientific advice meetings when QSP model is the main MIDD tool informing decision making. The technical focus of the MIDD meeting provides an optimal environment for in depth communication/discussion with the regulators on the CoU, model structure, assumptions, parameterization and validation. A detailed MAR (complete MAR) and a table for components for assessment of MIDD evidence (see complete table recommendation from ICH M15 (1)) are often included if an interim QSP analysis is performed. The Working Group advocates for the early adoption of the ICH M15 "Assessment of MIDD Evidence" table as a powerful tool to facilitate communication and engagement with both regulatory agencies and internal stakeholders. As a versatile methodology with an expanding model scope and CoU, QSP progresses alongside a drug's clinical development journey. Accordingly, the assessment framework for MIDD evidence should remain dynamic, adapting as data and knowledge accumulate and as the CoU evolves over time.

Model code, scripts and datasets for reproducing the modeling results were usually provided to the regulator as part of the MIDD meeting briefing package. The MIDD/scientific advice meetings can be used to share the proposed QSP modeling strategy that will be applied for the CoU, or it could be used to share interim QSP analysis to evaluate feasibility of a novel QSP approach to a CoU.

Sharing a draft version of model analysis plan (MAP) as part of the pre-filing activities may facilitate early communication with regulators and can provide clarity on how modeling results are planned to be used in cases where the CoU is novel for a QSP model. This is specially recommended if there is no peer review publication of the QSP model in the public domain that can be referenced. In the Working Group's past experiences, MAPs have been shared as appendixes to briefing books for type C meetings or MIDD meetings where the proposed QSP analysis was being discussed. The MAP should include a clear description of the CoU outlining how the model will support regulatory decision-making. It should also provide a succinct model description, including its biological scope, the translation of biological mechanisms into model structure, and the key underlying assumptions and limitations.

### 3.1.3 Reporting elements from case studies from submissions

The Working Group surveyed the FDA multidisciplinary review for NDA/BLA submissions containing QSP components. To complement these reviews (which do not disclose the full submitted reports) additional context was obtained from literature where it was possible to locate additional information on the models included in the multidisciplinary reviews. Six submission case studies were reviewed using rubric and the results are illustrated in Figure 2 and Appendix Figure 3. Although these submissions spanned different years, with some more current than others, all six reviewed collectively represent diverse therapeutic areas and modalities, summarized in Figure 2. In general, the associated reporting in those submissions were more comprehensive than those early-development QSP examples containing full model descriptions, limitations and schematics. While the model assessment table was not always initially included, they were frequently requested by FDA reviewers and subsequently provided, highlighting the importance of this assessment element when preparing QSP reports for submissions. Similarly, model codes were sometimes not initially provided with the original submission packages, but were shared upon request, indicating the importance of having annotated codes that can be easily shared with regulators and set up in a way that can be easy to review and run using typical software (Matlab scripts in many cases).

**1. Model Purpose & Regulatory Context**

| | Xenpozyme | Adzynma | Yorvipath | Natpara | Paxlovid | Elrexfio |
|---|---|---|---|---|---|---|
| Reg. Milestone | BLA | BLA | NDA | BLA | NDA | BLA |
| CoU | **Pediatric- adult extrapolation**<br>- Mechanistic similarity<br>- Efficacy in severe disease | **Supplement clinical data**<br>- Ultra rare disease and limited clinical data<br>- Build mechanistic links for disease relevant biomarkers | **Mechanistic E-R model**<br>- Explain observed E-R<br>- Evaluate untested dosing regimen. | **Provide dose justification**<br>- Compare dosing regimens<br>- Balancing efficacy and safety | **Support dose optimization**<br>- Justify treatment duration change in special population | **Support dose optimization**<br>- Evaluate impacts of sBCMA and dosing interval change<br>- Justify dose selection |
| Therapeutic Area | Rare/Metabolic | Rare | Metabolic | Rare/Metabolic | infectious diseases | Oncology (Heme) |
| Modality | Recombinant protein | Recombinant protein | Sustained-release prodrug | Peptide | Small molecule | Bispecific Antibody |

**2. Model Structure and Complexity**

| | Xenpozyme | Adzynma | Yorvipath | Natpara | Paxlovid | Elrexfio |
|---|---|---|---|---|---|---|
| Details & publications | Model details listed and in prior publications | Model details listed and in publications | Model details listed and in prior publications | Model details listed and in prior publications | Model details listed and in prior publications | Model details listed and in prior publications |
| Model Type | ODE model of drug MoA, mPBPK, disease pathophysiology, metabolism with multiple compartments | ODE model of drug MoA, disease pathophysiology with multiple compartments | ODE model of drug MoA, disease pathophysiology with multiple organ compartments. | ODE model of drug MoA, disease pathophysiology with multiple organ compartments. | ODE model of drug MoA, disease pathophysiology, immune response with multiple compartments | ODE model of drug MoA, disease pathophysiology, immune response with multiple compartments. |
| Scale | Molecular + Cellular + Tissue | Molecular + whole body | Molecular + Tissue + Organ | Molecular + Tissue + Organ | Cellular + Tissue | Molecular + Cellular |
| Assumptions & limitations | Included | Included | Included | Included | Included | Included |
| Model schematic | Multiple included | Multiple included | Multiple Included | Multiple Included | Multiple Included | Multiple Included |

**3. Model Development and Validation Report Components**

| | Xenpozyme | Adzynma | Yorvipath | Natpara | Paxlovid | Elrexfio |
|---|---|---|---|---|---|---|
| Data | In vitro, in vivo, natural history, clinical | In vitro, in vivo, clinical | Clinical data to guide literature model expansion | Clinical data to guide literature model modification | In vitro, in-vivo, clinical | In vitro, in-vivo, clinical |
| Codes | Model codes submitted | Could not be located | Model codes submitted | Published model code used | Model codes submitted | Model codes submitted |
| Calibration strategy | Staged calibration with publicly available data and internal individual clinical data | Staged calibration with internal clinical data | Staged calibration with publicly available clinical data and internal individual level clinical data | Staged calibration with subset of parameters fixed from literature & parameter input derived from clinical studies | Staged calibration w subset of parameters fixed from internal & literature data, paired with multiple Vpops | Staged calibration w subset of parameters fixed from internal & literature data, paired with multiple Vpops |
| Validation strategy | External healthy clinical data | External SOC datasets | Withheld clinical data; no out-of-sample predictive checks | Withheld internal clinical data | Withheld internal clinical data | Published external & emerging clinical data |
| Vpop/Twins | Vpop + Digital Twin | Digital Twins | No Vpop | No Vpop | Vpop | Vpop |
| UQ & SA | LSA & GSA & Vpops | LSA | Could not be located | LSA | LSA & Vpops | GSA & Vpops |

FIGURE 2 SUMMARY OF EVALUATION FOR SIX FDA MULTIDISCIPLINARY DISCLOSED REVIEWS

Overall, since at the BLA/NDA submission stage QSP models tended to be more mature, they had some version already published in the public domain, however details were provided in reports to describe the model in full and the calibration strategy in most reviewed cases. Sensitivity analyses were consistently included. Virtual populations (Vpops) or Digital twins were used in most cases, with only two exemptions (based on the available information from the multidisciplinary reviews). Overall, these elements reinforce the importance of a model uncertainty quantification strategy and the role that Vpops can have in model assessment and the need for inclusion in reports that support submissions. Notably, all evaluated cases had included staged calibration to various datasets both public and in-house as well as validation strategies.

The Working Group found the FDA multidisciplinary reviews and other publications from the agency highly valuable in informing the community of components of QSP submissions that seem to support MIDD evidence assessment as well as the nature of the feedback provided during review. The working group could not identify similar public resources with other regulatory authorities. It is our recommendation that global health authority organizations increase transparency and provide publicly components of QSP reports that have been received and how these submissions were evaluated given the CoU. Enhanced public disclosures of QSP submission will support alignment between regulators and the broader QSP community and advance best practices and expectations.

# 4 Proposed best practices for QSP model documentation and reporting

## 4.1 Sustainable QSP model documentation recommendations: internal technical peer review, code archiving practices

Best practices for building internal documentation components for MIDD rely on traceability, and iterative review early on rather than treating reporting as a final, stand-alone step. An important aspect is the use of internal review process at predefined development milestones, where model CoU, assumptions and limitations, are discussed, and aligned across internal stakeholders. Importantly, these discussions should be paired with formal QC of model code and analyses, ensuring consistency between scientific intent and technical implementation. A formal QC process should be built into the model lifecycle, not deferred to the end, and ideally conducted by qualified QSP modelers who were not directly involved in model development.

Robust code management is equally critical. Any project specific scripts should be version controlled and archived in regulatory compliant repositories, with clear documentation. Archiving should capture not only final code but also tagged milestone versions linked to internal decisions. Consistent file naming conventions, standardized folder structures, and readme files that explain how to reproduce results are good practice. Slide decks used in internal reviews should be archived alongside code to preserve model decision rationale and model history. With the adoption of the ICH M15 MIDD assessment table, QSP teams should consider using it not only at the end for reporting, when preparing to interact with health authorities, but as a planning, QC and documentation tool throughout development to facilitate alignment with internal stakeholders. Ultimately, high-quality reporting is only possible if the underlying scientific, and technical steps are documented systematically throughout the lifecycle of the QSP model and its intended CoU. As the QSP model's CoU changes as the drug development program advances, this comprehensive documentation of past CoU can be readily adapted to support formal regulatory communications.

## 4.2 QSP Model Analysis Report Components

To ensure transparency and efficient regulatory review, it is recommended that documentation of a QSP analysis in a dedicated report format is built gradually over time. Below is an outline of specific considerations for reporting components for QSP to be included in a MAR. Decisions regarding the inclusion and level of detail of each reporting element should be made by the QSP teams considering the novelty of CoU, type of regulatory interaction, and program stage. Other considerations should be aligned with level of model risk and impact as defined by the MIDD evidence assessment, and the publication history of the QSP model. The considerations presented reflect the Working Group's perspective on best practices for these specific reporting elements.

### 4.2.1 Introduction and objectives

In a QSP MAR, the introduction highlights the drug's targeted biological mechanism and its role in the relevant pathophysiology, outlining the biological processes represented, the rationale for the chosen model scope, and the constraints imposed by the CoU or available data.

*Model Purpose and Scope:* Most QSP MARs begin with a clear definition of the model's purpose and scope, outlining the biological and pharmacological processes it is intended to represent and specifying its regulatory QoI, such as dose selection, or trial design. Providing a model diagram is highly beneficial, as it offers a visual overview of the model structure including the key biological processes represented and helps reviewers quickly understand the biological components and their connectivity. Effective diagrams prioritize clarity by using consistent visual conventions (e.g., color coding, directional arrows, modular organization), minimizing unnecessary detail, and aligning naming conventions with those used in the model scripts.

*Mechanism of Action and Pharmacology Backbone*: Evaluation of a QSP model requires clear and accurate explanation of biological and pharmacological mechanisms included in the model. This should be supported by a concise summary of the current state of scientific literature, reflecting the totality of evidence (both in support of model assumptions as well as evidence that may not be as clear). Use of large language models (LLMs) can greatly facilitate this synthesis during model development; however, it is important that non-modeling subject matter experts SMEs (biology, clinical and clinical pharmacology colleagues) are engaged early for quality control of these components and review the MAR to ensure well rounded representation of totality of evidence. Such multidisciplinary SME support can help identify gaps in data or knowledge that may need to be handled through reasonable assumptions.

### 4.2.2 Data and Methods

*Data Sources: A* data section in a QSP MAR could be structured into two main elements for reporting: i) Available Data, ii) Dataset Preparation and Storage. For Available Data, it is recommended to clearly describe the data sources including clinical studies, and additional data from real-world evidence or literature. Relevant non-clinical data, such as preclinical in vivo studies, in vitro mechanistic experiments should also be listed.

For clinical studies, it is helpful to include a table of all the clinical trials, and natural history studies used to inform the QSP model, specifying key details like study phase, treatment information, biomarkers and endpoints used in the QSP analysis, participant numbers and any exclusion criteria. Where applicable, non-clinical studies that play a critical role in defining mechanistic relationships or parameter constraints may also be summarized in tabular or schematic form. In all cases, it is good practice to clearly indicate how each dataset was used whether for model development, calibration, validation, or qualitative support of mechanistic assumptions. Given the extensive data sources that are derived from published studies, it may not be feasible to describe in detail all the publications. It may be more practical to point to the sections in the MAR where these are described within the context of the model (such as parameter sources, model development, justification of assumptions).

For Dataset Preparation, this section could outline the software used for data analysis and detail the process of constructing analysis datasets. Specific procedures for handling missing data, samples below LLOQ (Lower Limit of Quantification), and outlier screening methods are helpful to include. Defining criteria for data exclusion and explaining how partial or missing data will be managed is also recommended.

*Model and Simulation Details:*

Model Assumptions and Limitations A dedicated assumptions section is recommended, with additional details placed in appendices where needed. All key assumptions should be explicitly stated and justified based on available biological, physiological, and pharmacological evidence, or clearly explained when they arise from data gaps or mathematical simplifications. Limitations should also be transparently described, including uncertainties, knowledge gaps, and validation constraints, to support appropriate interpretation of the model's reliability and boundaries.

Mathematical and Algorithmic Formulations: The mathematical formulation of the model (i.e. model variables, equations, fluxes, compartments, units etc) and any algorithms that are developed in a fit for purpose fashion to fit a model to datasets should be described clearly and stated in an easy to follow and detailed way either in a full technical MAR or an appendix. It is typically recommended that variable names are listed in separate tables and be chosen to be simple and easy to connect to any computational scripts. For this portion, it is our recommendation to treat a detailed MAR as a "modeling protocol" that is written to permit reproduction of model results. Since this part can be quite long, dedicated appendixes might help streamline MAR. While modeling scripts many times provide these types of details, it is important to include relevant mathematical derivations separately as to not obscure details due to syntax variation from one software to the next.

Model Parameters Parameter values are, in a way, the backbone of a QSP model. It is important that the values and their sources are summarized in a way that is easy to understand based on the model's CoU. It is common to summarize all this information in the form of a table, mostly in an appendix. A comprehensive table that includes parameter values, units, data sources, and their biological significance is essential for robust model evaluation. QSP parameters may be derived from diverse data types (in vitro, in vivo, ex vivo, clinical, or inferred) and their provenance must be considered when assessing their relevance and plausibility for representing the underlying biology. With these table formats, reliability of the sources is rarely evident, and evaluating each reference individually can be an enormous undertaking for models that have hundreds of parameters. The more reliable the sources of parameter values are, the more trust goes in the resulting simulations. For instance, receptor expression levels measured in immortalized cell lines may differ substantially from those obtained from patient-derived cells, with direct implications for defining biologically meaningful parameter ranges. Novel visualization methods have been introduced with a goal to provide a bird's eye view of the reliability and importance of sources of initial parameter values in a QSP model along with their influence in the model outputs. Health authority reviews of QSP analyses consistently emphasize the need for a rigorous assessment of parameter sources. Some of the case studies addressed this by applying these novel visualization methods to evaluate parameter reliability as described in Braakman et al (10).

In addition to documenting parameter values and sources, regulatory agencies increasingly require uncertainty quantification (UQ) of model parameters and their impact on predictions. It is important to distinguish uncertainty in a calibrated parameter estimate from population-level biological variability: the former reflects limited knowledge about the parameter value, whereas the latter reflects true between-subject heterogeneity and is often more relevant for clinical interpretation. For parameters estimated through calibration, uncertainty may be reported as confidence intervals or coefficient of variation (CV%), while for fixed parameters derived from literature or experiments, biologically plausible ranges and their rationale should be documented.

The combined impact of parameter uncertainty and biological variability on model predictions can be evaluated using virtual populations (Vpops) that are calibrated to match the observed data and population variability, thereby providing more constrained and credible predictions than simple ensemble simulations based on independent sampling of a single parameter. This approach can provide prediction intervals rather than point estimates. The resulting prediction intervals (e.g., 90% or 95% intervals) communicate the range of plausible outcomes given parameter uncertainty and population-level biological variability and should be reported for all key model outputs relevant to the regulatory question.

Model Simulation: This section should clearly describe the simulations conducted with the model, including the rationale for selecting specific scenarios and how each simulation supports the stated CoU. For each scenario, the underlying assumptions, inputs, and intended questions should be explicitly articulated to ensure transparency and interpretability of the results. When context-specific parameter sets are employed, (for a specific scenario or Vpops) their definition, derivation, and scope of applicability should be documented explicitly. This includes clarifying whether parameter values are shared across scenarios or tailored to reflect specific populations, disease or severity states, or treatment conditions. Such documentation enables reviewers to understand how simulation results are generated and how they support the defined CoU.

Model Workflow: Given the often complex and iterative nature of QSP model development, such diagrams would be particularly helpful during regulatory review. For instance, diagrams describing model calibration and validation strategies, the sequencing and integration of different datasets used to inform these activities, and the criteria for model refinement. In addition, visual representations of virtual population generation and qualification strategies, such as parameter sampling approach, acceptance criteria, and links to intended CoU are strongly recommended to enhance interpretability during review.

### *Approaches for Model Evaluation:*

Model Verification: Verification in QSP models involves several complementary activities aimed at ensuring the correct implementation of the model for the CoU. Code verification focuses on checking that equations, units, and parameter definitions are implemented correctly. In parallel, model structural checks should be done to evaluate whether the biological processes represented mathematically behave consistently with the underlying biological assumptions, without any issues related to mass balance, or numerical instabilities.

Model Validation: Validation focuses on whether the model's predictions reliably represent biological and clinical behavior under the defined CoU. This is an important step that allows regulators to assess how reliable and appropriate the model's predictions are for the decision being informed. Thus, validation strategy should be context-driven and proportionate to model influence as described in the MIDD evidence assessment, ensuring that the rigor applied aligns with the potential impact of model-informed decisions (11). Validation in QSP models often frequently encompasses multiple biomarkers, pathway dynamics, clinical endpoints, and is therefore commonly conducted using VPop-based comparisons to evaluate model performance considering variability rather than a single "typical" trajectory. VPop construction is essentially an inverse calibration process in which parameter distributions are identified to produce simulated outcomes matching observed clinical responses. Because VPops capture inter-individual variability, documenting their generation is essential for clarity and reproducibility. Key reporting elements include the parameter sampling strategy used, parameter ranges and biological justification, acceptance criteria for virtual patients, and any refinement steps such as filtering or reweighing to better align with observed population characteristics, along with metrics used to assess that alignment. Validation typically includes several interrelated components that need to be reported, such as the validation criteria, selecting and comparing against appropriate datasets, and conducting both qualitative and quantitative evaluation of model performance.

Validation Criteria Clear, predefined validation criteria are essential and should align directly with the model's CoU and decision risk. Acceptability may be based on prediction error margins or agreement with observed variability, with a clear rationale for why each criterion is appropriate given model risk. When model-based decisions carry higher consequences, agencies may expect validation using clinical data from the target population.

Appropriate validation datasets External validation, testing the model against data not used during model development, is important to establishing credibility, as it demonstrates predictive capability beyond the calibration data. In regulatory contexts, demonstrating model performance on independent datasets, particularly those relevant to the intended population and CoU, significantly strengthens the case for its reliability (11).

Performance evaluation The evaluation of validation outcomes typically includes a combination of graphical and statistical diagnostics that characterize the model's predictive adequacy and residual uncertainty. Visual Predictive Checks (VPCs) are a standard component of model performance evaluation and can be effectively applied to QSP -derived Vpops. Agreement between observed outcomes and model-predicted intervals supports the model's ability to reproduce population-level trends and variability. See Figure 3 for examples from published FDA multi-disciplinary reviews. Traditional GoF diagnostics, such as observed versus predicted values, residual plots, or bias assessments, also provide valuable information about model adequacy, though often at a more aggregate level in QSP models, Figure 3.

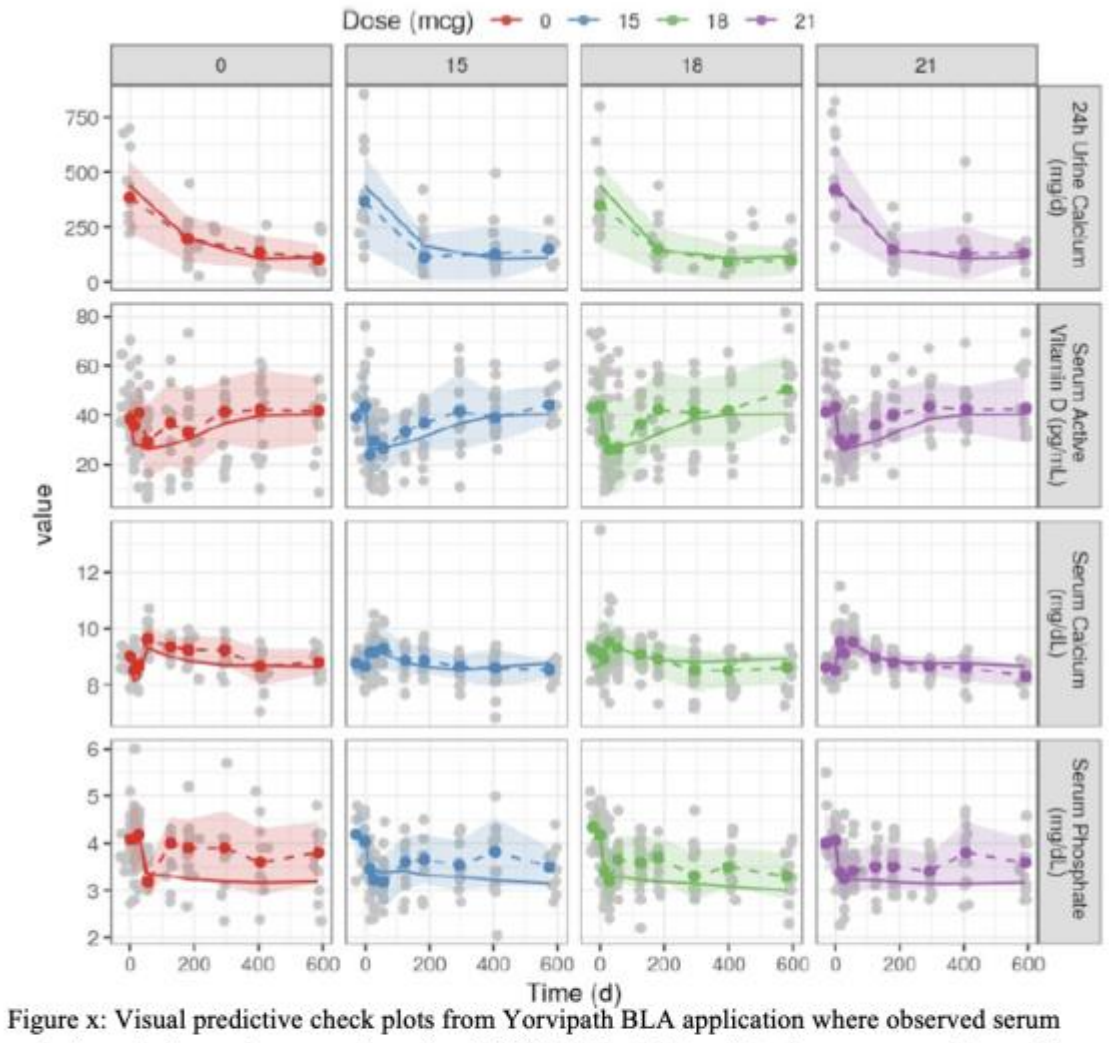

Figure x: Visual predictive check plots from Yorvipath BLA application where observed serum and urine calcium and serum active vitamin D (individual data points shown as gray dots, median as colored dashed line) were compared to model predictions (solid lines)

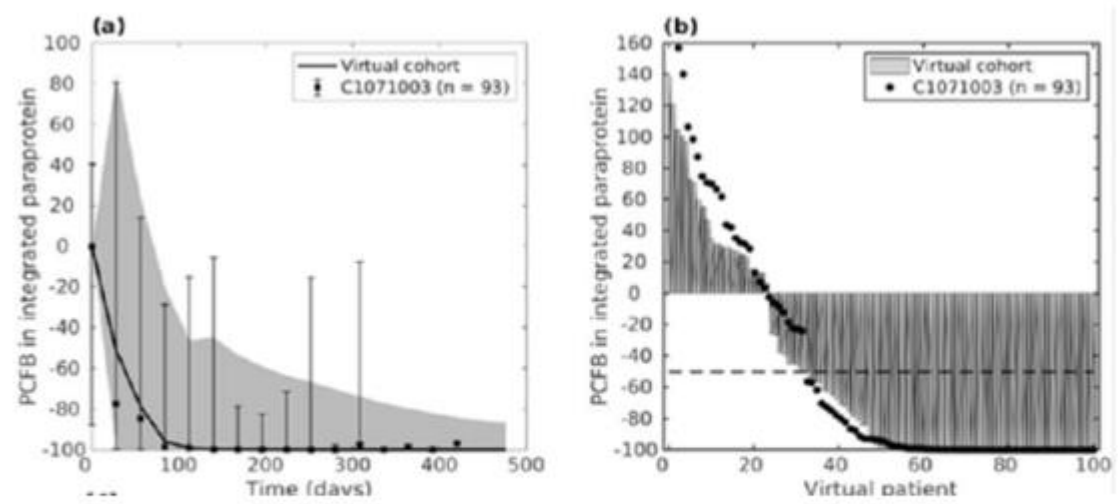

Visual predictive check plots from Elrexfio BLA application that compares virtual cohort simulations with C1071003 clinical trial results. (a) Spider plot illustrating percentage change from baseline (PCFB) over time. The solid black curve depicts the median trajectory across the virtual population (VPop), while the grey shaded region represents the corresponding 95% prediction interval. Observed pooled participant data are summarized by black square markers indicating the median, with vertical error bars showing ±2 standard errors.

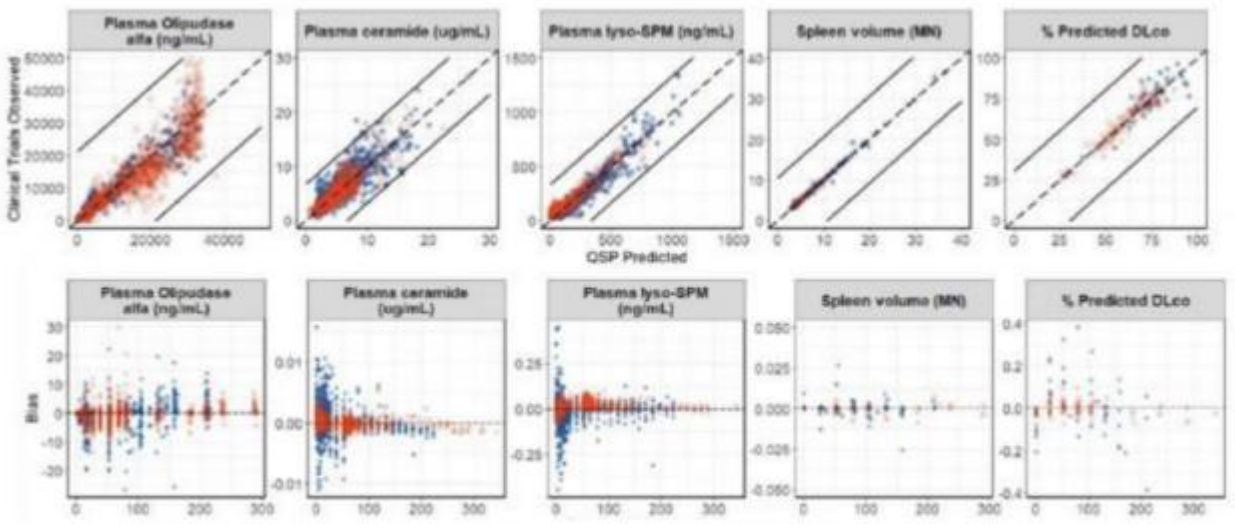

Goodness of Fit and Bia Plots Comparing Predicted Model Outputs to the Corresponding Observed Data in XENPOZYME BLA application. Black solid lines are 2σ of clinical data from all patients included in QSP analysis. Blue dots represent pediatrics, red dots represent adults

**Figure 3Examples of graphical evaluations of validation of QSP models**

## 4.2.3 Results

An important component of QSP model documentation in regulatory filings involves communicating simulation results in a way that is informative for the question of interest. Ideally, reviewers can follow what was simulated, why particular scenarios were explored, and how these outcomes relate to the broader drug development program. Several elements often contribute to achieving this clarity, including articulation of baseline and treatment scenarios, characterization of the virtual population used in the analysis, and the integration of visualization and discussion to convey both simulation findings and their clinical relevance.

**Table 4 Simulation Scenarios and their reporting elements**

| | Purpose | Key Elements to Report | Regulatory Relevance |
|---|---|---|---|
| Baseline (Untreated) Scenario | Characterize the natural history of disease in the absence of therapy. | • Virtual placebo arm or simulated disease progression curve<br>• Assumptions regarding baseline disease severity and biomarker dynamics<br>• Time-course trajectories (e.g., median with prediction intervals) clearly labeled as untreated | Anchors interpretation of treatment effects and allows assessment of physiological plausibility, particularly in chronic or progressive diseases where untreated biomarker trends provide essential context for predicted treatment effects. |
| Treatment Scenarios | Evaluate drug effects under defined clinical assumptions. | • Dosing regimen (dose, frequency, duration)<br>• Baseline population characteristics<br>• Concomitant therapies or combinations<br>• Clear distinction between simulations used for model alignment/calibration and those exploring new doses or populations | Clarifies the evidentiary basis of simulations and supports interpretation of whether results are confirmatory (aligned to existing data) or extrapolative (prospective predictions). |
| Sensitivity / Alternative Scenarios | Assess robustness of conclusions and explore uncertainty or mechanistic alternatives. | • Parameter uncertainty analyses<br>• Simulations under extreme physiological conditions<br>• Alternative mechanistic assumptions<br>• Clear annotation of "what-if" purpose and limitations | Supports MIDD evidence assessment and contextualizes uncertainty; ensures reviewers understand the scope and limitations of exploratory simulations. |
| Virtual Population Characteristics | Describe the population framework used for all simulations. | • Demographic and physiological assumptions<br>• Baseline disease characteristics<br>• Ranges and distributions of key virtual parameters (summarized in tables or visualized via box/violin plots)<br>• Relevance of VPop assumptions to the scientific question | Ensures transparency regarding population representativeness and enables assessment of clinical relevance and generalizability of simulation outcomes. |

### 4.2.4 Discussion and conclusions

A discussion of clinical relevance often complements the presentation of results by articulating how the simulations contribute to understanding the therapeutic profile of the drug candidate. Such a discussion may draw connections between the model predictions and anticipated clinical outcomes, noting how the simulations inform considerations such as dose justification, potential safety concerns, or expected responses in populations for which empirical data may be limited. Because QSP models can explore dosing strategies, subpopulations, or mechanistic hypotheses that are difficult to evaluate experimentally, this discussion provides an opportunity to reflect on how the simulations extend or contextualize available clinical evidence.

It can be valuable to acknowledge uncertainties or discrepancies between simulated and observed outcomes, considering possible mechanistic explanations and whether these have implications for the model's suitability within its specified CoU. The section may conclude by summarizing how the collective findings address the regulatory question of interest and by noting any limitations that may influence the confidence placed in the predictions. Framing the discussion in this manner helps ensure that simulation results are not only reported but meaningfully connected to the clinical and regulatory decisions they are intended to inform.

# 5 Conclusions and Future Directions

Establishing consistent guidance and standardized reporting elements is critical for ensuring QSP model credibility, which is essential for regulatory confidence and scientific rigor. This white paper proposes best practices in reporting based on QSP analyses that have impacted successful regulatory approvals. These best practices also have considered early regulatory reporting considerations based on de-identified programs, along with the experience of WG members. All this information was grounded on the available reporting recommendations derived from PBPK guidelines and the recent ICH M15 guidance.

Based on these, a tiered approach to QSP reporting is recommended as shown in Figure 4. This approach defines a set of minimal reporting elements, with the inclusion of additional elements determined by factors such as the stage of development, the novelty of the CoU, and the level of model risk. At a minimum, QSP MAR should include an introduction and objectives, assumptions and limitations, data sources, a description of model development and methods, simulation results, and conclusions. The depth of reporting can be further adjusted based on the model's peer-reviewed publication history, with models supported by more extensive publication histories requiring a more streamlined MAR.

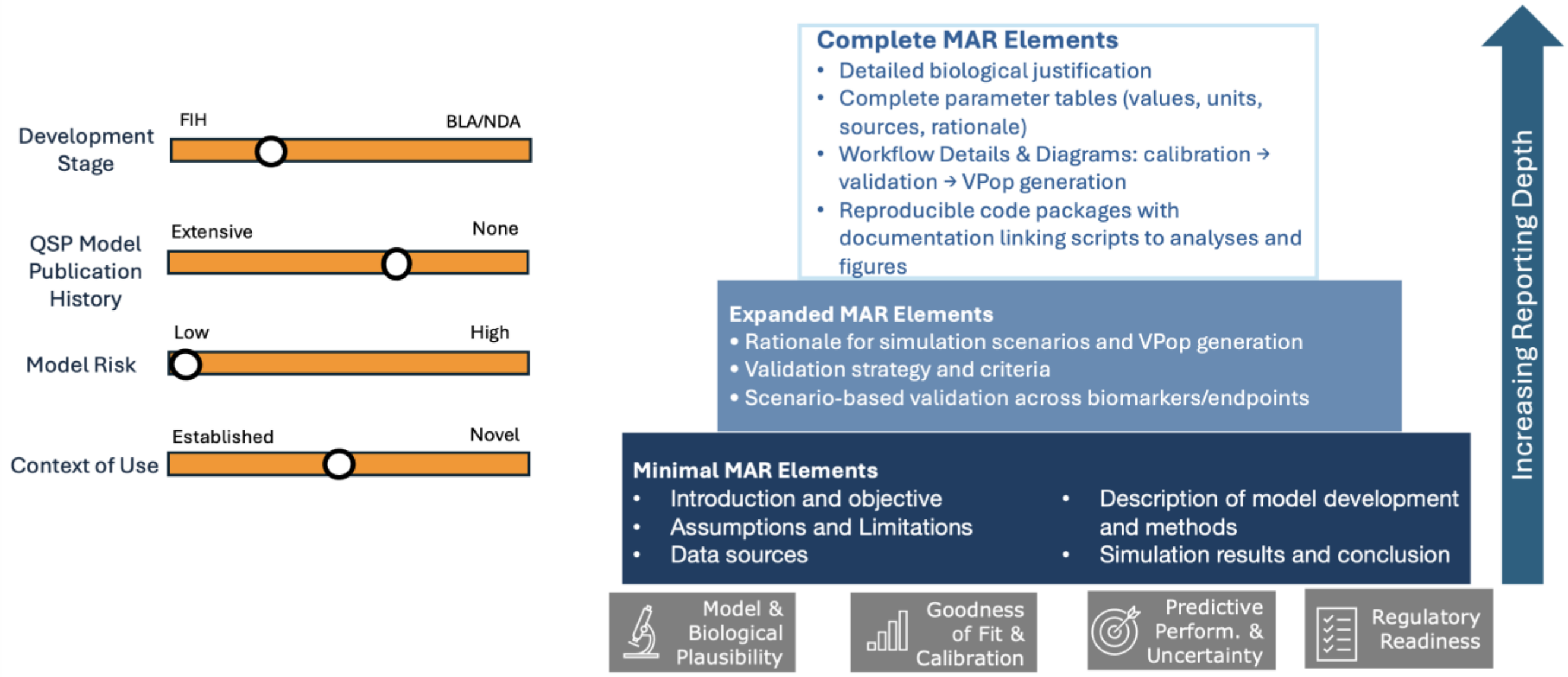

FIGURE 4 SUGGESTED TIERED APPROACH FOR QSP REPORTING

Reporting of QSP models is the foundation to ensuring reproducibility and transparency, which are central to establishing regulatory confidence in MIDD-derived evidence. Clear, structured documentation enables regulators to understand the model's CoU, key assumptions, data sources, and uncertainties, and to assess whether the model is fit for its intended decision-making purpose. Embedding these principles within a reporting framework, consistent with ICH M15 guidance, is essential for advancing QSP from a supportive analytical tool to a routinely accepted and reliable component of drug development and regulatory decision-making.

Equally important is continued investment in education and training to build shared understanding among modelers, clinicians, and regulators of the value and limitations of QSP modeling. Strengthening the ability of stakeholders to appropriately interpret QSP results is critical for effective use of these models in drug development and regulatory decision-making.

While these best reporting practices build on lessons learned from previous QSP applications, it is critical to adopt a forward-looking perspective that anticipates evolving regulatory expectations. Future credibility frameworks should not only reflect historical successes but also proactively address how regulators will assess QSP models in the near term, particularly as AI-driven tools could enable more efficient MAR generation and review processes.

# 6 Conflict of Interest

BS is an employee of Pfizer Inc. and may hold stock or stock options. SZ is an employee of Sanofi and may hold Sanofi stock or stock options. YC is an employee of Daiichi Sankyo and may hold Daiichi Sankyo stock or stock options. AG is employed by and may own stocks in Merck & Co., Inc. WW is an is an employee of Johnson & Johnson and may own shares or stock options in Johnson & Johnson. AVR is an employee of Bristol Myers Squibb and may hold stock or stock options. FAS is an employee of Genmab A/S and may hold stock or stock options. JQXG is an employee of GSK and may hold stock or stock options in GSK. AK is an employee of ESQlabs GmbH. MES is an employee of Genentech Inc and may hold stock or stock options. CS is employed by and may own stocks or stock options in Simulations Plus Inc. The remaining authors declare no competing financial interests.

# 8 Appendix

| QSP Case Studies Framework | Suggested Details |
|---|---|
| **1. Model Purpose/Regulatory Context** | |
| Primary objectives/CoUs | |
| Regulatory milestone | IND, NDA, BLA, etc. |
| Interaction/impact type | A/D, EoP2 meetings |
| Therapeutic area and indication(s) | |
| Drug Modality | |
| **2. Model Structure and Complexity** | |
| Level of detail if Published | Cited/described? |
| Model type | Mechanistic depth (semi-mechanistic vs. comprehensive platform model) |
| Biological scale | Molecular, cellular, tissue, organ, whole-body integration |
| Assumptions and limitations | Listed? |
| Model schematic | Y/N |
| **3. Model Development and Validation** | |
| **Data sources** | Types of data incorporated (in vitro, preclinical, clinical) |
| Code verification | |
| Calibration strategy | Parameter estimation approach and calibration datasets |
| Validation strategy | Internal vs. external validation, prospective validation, and metrics |
| Vpop/digital twin generation | Approach detailed if internal example? Metrics for qualifying VPop? |
| Uncertainty quantification and parameter sensitivity analysis | |
| **4. Results and Conclusions** | |
| Were model results aligned with Context of Use and key questions | |
| Were model limitations contextualized with results | |
| **5. Regulatory Feedback and Impact** | |
| Credibility Assessment | Alignment (sponsor vs reviewer) |
| Regulatory feedback | Nature of regulatory questions and concerns |
| Reproducibility of simulations | Further analysis by reviewer |
| Regulatory acceptance | Level of regulatory buy-in for model-based conclusions |

APPENDIX FIGURE 1 CASE STUDY EVALUATION RUBRIC

| 1. Model Purpose/Regulatory Context | | |
|---|---|---|
| **Regulatory Milestone & Interaction Type** | IND, minimal MAR | Ph2b Briefing Book |
| **CoU** | FIH dose regimen | Combination Study Design |
| **Therapeutic Area & Modality** | Immunology, mAbs | Immunology, mAbs |

| 2. Model Structure and Complexity Report Components | | |
|---|---|---|
| **Details & Publications** | Model not disclosed at the time of IND. Minimal MAR was included that detailed biology included in model, assumptions and references. | Model components, key assumptions listed<br>An abbreviated version of the model planned for publication |
| **Model type** | QSP model described drug MoA and disease pathophysiology, connecting to biomarkers and clinical endpoint. | Full-scale, disease-centric QSP model that incorporates all key cytokines, cell types known to contribute to disease pathology as well as clinical endpoints |
| **Biological scale** | Subcellular + Cellular + Tissue | Cellular + Tissue |
| **Assumptions and limitations** | High level assumptions and limitations included in report | Assumptions and limitations were captured in report |
| **Model schematic** | Single | Multiple included |

| 3. Model Development and Validation Report Components | | |
|---|---|---|
| **Data sources** | In vitro, in-vivo, natural history data, clinical data from published studies from approved drugs in same indication included | In vitro, pre-clinical and clinical data; clinical data encompass all major therapies across multiple MoAs, including both in-house individual-level patient data and published aggregated-level data |
| **Code verification** | Model codes not submitted | Model codes not submitted |
| **Calibration strategy** | Brief description of model calibration | Described in a report: parameter derived from literature, in vitro, non-clinical data. Virtual populations calibrated to clinical outcome data across multiple therapies. |
| **Validation strategy** | Validation against data from approved drugs in same indication. | For each key MoA, clinical datasets were pre-divided into training and validation subsets. A priori prediction of novel compounds and novel mechanisms also performed as part of the model assessment. |
| **Vpop/digital twin** | Vpop | Vpop |
| **UQ & SA** | LSA | LSA and Vpop |

| 4. Results and Conclusions Report Components | | |
|---|---|---|
| **Model results aligned with CoU & key questions** | Model justified the recommended dose, regimen for FIH trial | Model supported dose selection and combinations for a Phase 2 POC study |
| **Model limitations contextualized** | Model limitations stated in report | Model limitations stated in the report |

| 5. Regulatory Feedback and Impact Report Components | | |
|---|---|---|
| **Risk Assessment** | Sponsor and review team agreed on: Model Impact : **Med**.<br>Decision Consequence: **Low**. | Sponsor proposed: Model Impact : **High**.<br>Decision Consequence: **Low**.<br>No objections were received from the review team |
| **Regulatory feedback** | NA | Requested to provide more detailed modeling report in the future |
| **Reproducibility of simulations** | NA | No additional model analysis and codes was requested |
| **Regulatory acceptance** | NA | NA |

APPENDIX FIGURE 2 ANONYMIZED CASE STUDIES FROM WORKING GROUP

**4. Results and Conclusions Report Components**

| | Xenpozyme | Adzynma | Yorvipath | Natpara | Elrexfio | Paxlovid |
|---|---|---|---|---|---|---|
| Were model results aligned with Context of Use and key questions | QSP analysis supported similarity of disease between pediatric and adult patients, justifying extrapolation of efficacy | Yes. Model assessed the mechanistic relationships among ADAMTS13/VWF activities and platelet count to supplement limited clinical data | QSP model supported mechanistic understanding of the E-R relationship. | Model simulations supported testing a more frequent than QD or slow-release QD dosing to provide better control of hypercalciuria | Model supported recommended dose and regimen registered in the approved label for high and low baseline soluble BCMA RRMM patients | Model justified the recommended dose, regimen and treatment duration for clinical trials in EPIC-HR and immunocompromised patients. |
| Were model limitations contextualized with results | Model limitations stated in publications, report and discussed during review | Model limitations thoroughly discussed in the publication | Model limitations stated in the report and discussed during review | Model limitations stated in publications, report and discussed during review | Model limitations stated in publications, report and discussed during review | Model limitations stated in publications, report and discussed during review |

**5. Regulatory Feedback and Impact Report Components**

| | | | | | | |
|---|---|---|---|---|---|---|
| Risk Assessment | Sponsor and review team agreed on: Model Impact : Low. Decision Consequence: Low. | The BLA review suggest possible Moderate for both Model Impact and Decision Consequence, though they were not explicitly mentioned. | Sponsor and review team agreed on: Decision Consequence: Low. Model Influence: Sponsor proposed as Low, review team assessed as Medium as QSP modeling is considered as an alternative to E-R analysis | Not performed (Review in year 2013) | Sponsor and review team agreed on: Model Impact : Low. Decision Consequence: Low. | Sponsor and review team agreed on: Model Impact : Medium. Decision Consequence: Low. |
| Regulatory feedback | Requested additional specifications to facilitate independent model script runs. | A pharmacometrics QSP consult review was conducted, and it agreed that the model provided confirmative evidence. | Requested code and additional specifications to facilitate independent model script runs. Early communication of analysis documentation will help with timely review. | Not specified in the review | Requested all codes and datafiles, hardware definition files and specifications to facilitate independent model script runs | Requested all codes and datafiles, hardware definition files and specifications to facilitate independent model script runs |
| Reproducibility of simulations | Additional model analysis and scripts were requested by reviewers to evaluate model performance and assess sensitivity of model predictions to calibration strategy, as detailed in BLA review. | Not specifically mentioned in BLA review. | Additional model analysis and scripts were requested by reviewers to evaluate model performance. Reviewer can reproduce the results in the report. | Reviewer independently reconstructed the model and performed simulations to evaluate typical patient and alternative dose scenarios. The results supported | Additional model analysis and codes were requested by reviewers to evaluate model performance and assess sensitivity of model predictions to certain assumptions, as detailed in BLA review. | Additional model analysis and codes were requested by reviewers to evaluate model performance and assess sensitivity of model predictions to certain assumptions, as detailed in NDA review. |
| Regulatory acceptance | Model was accepted as supportive evidence both in its ability to capture clinical data of interest (Quality and robustness of QSP model and fitting strategy) and to be able to answer both questions in CoU. | Model was accepted. FDA explicitly states the model provides confirmative evidence to support the approval of ADZYNMA in cTTP. | Model was accepted as supportive evidence for the 1st CoU in its ability to capture clinical PK and biomarker data for disease aspects and treatment response The 2nd CoU was not accepted due to lack of out-of-sample validation for higher dose groups in patient model validation. | Model was accepted as supportive evidence for the CoU, supported evaluation of alternative dosing strategy, that warrant further clinical investigation. | Model was accepted as supportive evidence for both Context of Use cases having shown acceptable calibration for PK and PD components with provided clinical datasets and overall mechanism of action. | Model was accepted as supportive evidence both in its ability to capture clinical data of interest (Quality and robustness of QSP model and fitting strategy) and to be able to answer both questions in CoU. |

APPENDIX FIGURE 3 PUBLISHED FDA MULTIDISCIPLINARY REVIEW CASE STUDIES (CONTINUED)

## 8.1 BLA/NDA Case Study Sources

#### *8.1.1.1 Xenpozyme*

https://www.accessdata.fda.gov/drugsatfda_docs/nda/2022/761261Orig1s000IntegratedR.pdf

#### *8.1.1.2 Adzynma*

https://www.fda.gov/media/174318/download

https://www.ema.europa.eu/en/documents/assessment-report/adzynma-epar-public-assessment-report_en.pdf

McBride, C., Jiang, J., Zhang, Z., Tolsma, J., Patwari, P., Mellgård, B., Vakilynejad, M., Bhattacharya, I. and Zhu, A.Z.X. (2025), Quantitative Systems Pharmacology Modeling of Platelet Responses to Recombinant ADAMTS13 in Patients With Congenital Thrombotic Thrombocytopenic Purpura. CPT Pharmacometrics Syst Pharmacol, 14: 1575-1586.

#### *8.1.1.3 Elrexfio*

https://www.accessdata.fda.gov/drugsatfda_docs/nda/2023/761345Orig1s000MultidisciplineR.pdf

*8.1.1.4 Paxlovid*

https://www.accessdata.fda.gov/drugsatfda_docs/nda/2023/217188Orig1s000IntegratedR.pdf

*8.1.1.5 Natpara*

https://www.accessdata.fda.gov/drugsatfda_docs/nda/2015/125511Orig1s000ClinPharmR.pdf

*8.1.1.6 Yorvipath*

https://www.accessdata.fda.gov/drugsatfda_docs/nda/2024/216490Orig1s000MultidisciplineR.pdf